\documentclass[a4,12pt]{article}
\usepackage{bbm}
\usepackage{amsmath}
\usepackage{amssymb}
\usepackage{amsfonts}
\usepackage{cite}
\usepackage[dvips]{graphicx}

\numberwithin{equation}{section}

\newtheorem{theorem}{Theorem}[section]

\newtheorem{conjecture}[theorem]{Conjecture}
\newtheorem{corollary}[theorem]{Corollary}
\newtheorem{defn}[theorem]{Definition}

\newtheorem{lemma}[theorem]{Lemma}
\newtheorem{prop}[theorem]{Proposition}

\newtheorem{remark}[theorem]{Remark}

\def \mc{\mathcal}

\def \v{\vskip 0.1in}
\def \n{\noindent}

\def \real{\mathbb{R}}

\def \mf {\mathfrak}

\begin{document}

\title{Uniform K-stability for extremal metrics on toric varieties}

\author{Bohui Chen, An-Min Li, Li Sheng}

\maketitle

\abstract In this paper we prove that for toric varieties the uniform
K-stability is a necessary condition for the existence of extremal metrics. 
 \endabstract

\section{Introduction}
Extremal metrics, introduced by E. Calabi, have been
studied intensively in the past 20 years. The necessary conditions
for the existence are conjectured to be related to certain
stabilities. There are many works on this aspect. In particular,  Tian first made great progess toward understanding this.
It was   Tian who first  gave an analytic "stability"
condition which he proved that is equivalent to the existence of a K$\ddot{a}$hler-Einstein metric (\cite{T1}). This condition is
the {\em properness} of the Mabuchi functional.
In \cite{T1}, Tian also defined the algebro-geometric notion of K-stability.
In \cite{D1}, Donaldson generalized Tian's definition of K-stablilty by giving
an algebro-geometric definition of the Futaki invariant and conjectured that
it is equivalent to the existence of a cscK metric. {  The conjecture can be formulated as
\begin{conjecture}[Yau\cite{Y},Tian\cite{T1},Donaldson\cite{D1}]\label{conj_1}
The manifold $M$ admits a cscK metric in the class $c_1(L)$ if and only if $(M,L)$ is K-polystable.
\end{conjecture}
In \cite{T} and \cite{T2}, Tian explored the extremal metrics and geometric stability and conjectured
that $(M,\Omega)$ admits a cscK metric if and only if $(M,\Omega)$ is {\em analytically $G$-stable}
for some maximal compact subgroup $G$ of $Aut(M,\Omega)$ (cf. Conjecture 4.1(\cite{T2})). In \cite{St},
Stoppa has a similar conjecture that $(M,L)$ admits an extremal metric in the class $c_1(L)$ if and
only if the pair $(M,L)$ is K-polystable relative to a maximal torus of $Aut(M,L)$, while the necessary condition of this conjecture
 is
proved in \cite{SS}.}

The problem of searching canonical metrics may become  simpler if the manifold studied
admits more symmetry. Hence, it is natural to consider the  toric varieties first.
Each toric manifold $M^{2n}$ can be represented by a Delzant
polytope $\bar \Delta$ in $\real^n$, the
equation for the extremal metrics becomes a real 4th order equation which is known as the
Abreu equation(cf. \eqref{eqn_1.4}).
In a sequence of papers, Donaldson
initiated  a program to study the extremal metrics on toric manifolds.
The problem is to solve the equation under certain
necessary stability conditions.

The formulation of stability becomes more
elementary when we consider the toric varieties. We list some important results on this topic:
\begin{enumerate}
\item In \cite{D1},
Donaldson
formulated {\em relatively K-polystability} for Delzant polytopes and conjectured  that this
stability implies the existence of the cscK metric on toric manifolds (cf. Conjecture \ref{conj_1}).
\item  In \cite{D1}, Donaldson also introduce a
stronger version of stability which we call   {\em uniform \rm{K}-stability} in this paper(cf. Definition \ref{definition_1.2.4}). He observed that
 the existence of weak solutions follow immediately  from such a
uniform   stability(\cite{D1}, also Theorem \ref{theorem_4.7} in \S\ref{sect_4.3}).
\item
When n=2, relatively K-polystability  may imply uniform K-stability:
in \cite{D1}, Donaldson prove that the relative K-polystability imply the uniform K-stability when  $A>0$;
recently, in \cite{WZ}, Wang-Zhou show that  the relative K-polystability implies the uniform K-stability when $A$  is an affine linear function.
\item
Zhou-Zhu(\cite{ZZ1}) introduced the notion of {\em properness}, i.e, the analytic stability,  on the modified
Mabuchi functional and showed the existence of weak solutions under
this assumption.
\end{enumerate}
On the other hand, In \cite{ZZ2}, Zhou-Zhu proved that relative K-polystable
is a necessary condition.

In the thesis of Szekelyhidi(\cite{Sz}), the author gave a nice survey and exposition of K-stability.
In particular, one version of uniform stability is discussed (cf. \S3.1.1(\cite{Sz}), also the discussion at the
end of \S\ref{sect_4.2}) and correspondingly,
such a uniform stability on toric manifolds is studied in \S4.2.1(\cite{Sz}). As mentioned in his thesis,
"it  is tempting to conjecture that uniform K-stability is the correct condition characterizing the existence
of cscK metrics". Such a point reveals the importance of uniform stability.

In this paper,
our goal is to show that the uniform K-stability (introduced by Donaldson(\cite{D1})) {is} a {\em necessary}
 condition of the existence of extremal metrics
{\em on toric varieties} (cf. Theorem \ref{theorem_4.6}).
Such a uniform stability is defined in Definition \ref{definition_1.2.4}, it
is stronger than the one given in \cite{Sz}. As this  uniform K-stability immediately
implies the existence of weak solution of extremal metric, we would suggest a modification on the original conjecture posed by Donaldson:
\begin{conjecture}
Let $(M,\omega)$ be a toric manifold and $\Delta$ be the corresponding Delzant polytope.
$(M,\omega)$ has a metric within the class $[\omega]$ that solves the Abreu equation
\eqref{eqn_1.4} if and only if $(\Delta,A)$ is uniformly K-stable. In particular,
the metric is extremal if and only if $A$ is affine linear.
\end{conjecture}
Hence, our main theorem \ref{theorem_4.6} says that the direction "$\Rightarrow$" is true.
On the other hand, by Donaldson's result, uniform K-stability implies the existence of  the weak solution (see also Theorem \ref{theorem_4.7}).
In order to prove the conjecture, it remains to establish the regularity of the weak solutions. On the other hand,
as a consequence of uniform stability, we have the properness of the Mabuchi functional(cf. Theorem \ref{theorem_4.8}).
Hence, we show that the properness of the Mabuchi functional is a necessary condition for the extremal metrics on toric
varieties, which is part of Tian's conjecture on toric manifolds.
\begin{remark}
(i) By settling the regularity theory for toric surfaces, we are able to prove the conjecture
for almost any $A$, (\cite{CLS1} and \cite{CLS2}). (ii) It is tempting to
get $C^0$ estimate for the weak solution from the uniform K-polystability, this is done by
Donaldson when $n=2$. We will consider this in the further study. Note that, by the work of \cite{CLS3},
such an estimate would imply the interior regularity of the Abreu equation for any dimension.
\end{remark}

There is a natural generalization of this problem in the PDE sense. Namely, we
may consider $\Delta$ to be arbitrary convex domain and study the variational problem for the Abreu equation. The stability
arises naturally for such set-ups and it appears to be  the necessary condition for minimizer.   This is considered by
Donaldson as well and we give a short discussion on this generalization in \S\ref{sect_2}. On the other hand, there is a family
of 4th order PDEs which are of  the same type of the Abreu equation (cf. \cite{LJSX}). A generalization of the theory in this paper to those equations
will be discussed in \cite{CLS4}.

\section{A variational problem}\label{sect_1.1}

Let $\Delta$ be an open polytope in $\mf t^\ast\cong \real^n$. Suppose its closure $\bar\Delta$ is compact.
We explain the variation problem following  Donaldson's paper closely (cf. \S 3 in \cite{D1}).

Suppose that
$\Delta$ is defined by linear inequalities $h_k(x)-c_k>0$, where $h_k(x)-c_k=0$ defines the facet of $\Delta$. Write $\delta_k(x)=h_k(x)-c_k$
and set
\begin{equation}\label{eqn_1.0}
u_o(x)=\sum\delta_k(x)\log\delta_k(x).
\end{equation}
This function is first defined by Guillemin and it defines a K\"ahler metric on the toric variety defined by $\Delta$ if it is a Delzant polytope.

Define
$\mc S$ to be the set of continuous convex functions $u$ on $\bar\Delta$ such that $u-u_o$ is smooth on $\bar\Delta$.
Let $A$ be a {\em smooth} function on $\bar \Delta$ (in fact, we may only require that $A\in L^\infty(\Delta)$ (cf. \cite{D1}).
Consider the following functional
\begin{equation}\label{eqn_1.1}
\mc F_A(u)=-\int_\Delta \log\det(u_{ij})d\mu+\mc L_A(u)
\end{equation}
where
\begin{equation}\label{eqn_1.2}
\mc L_A(u)=\int_{\partial\Delta}ud\sigma-\int_\Delta Au d\mu.
\end{equation}
$\mc F_A$ is known to be the Mabuchi functional and $\mc L_A$ is closely related to the Futaki invariants.

\begin{defn}\label{defn_1.1}
Let $u$ be a continuous convex function on $\Delta$. We say it is an {\em extremal function}
if $\mc L_A(u)=0$.
\end{defn}
\begin{remark}\label{remark_1.1.0}
It is easy to see that the functional $\mc F_A(u)$ is bounded below only if that any affine linear function is
an extremal function. Hence,
We always assume that any affine linear function is an extremal function. This is equivalent to
the vanishing of the Futaki invariants. Let $\mc G$ denote the set of affine
linear functions.
\end{remark}
An extremal function is called  {\em nontrivial} if it is not affine linear, i.e, not in $\mc G$.
By this assumption, it is clear that the functional $\mc F_A$ reduces to be a functional on
$\mc S/\mc G$.

The Euler-Lagrangian equation for $\mc F_A$ is
\begin{equation}\label{eqn_1.4}
\sum\frac{\partial^2 u^{ij}}{\partial x_i\partial x_j}=-A.
\end{equation}
This is called the Abreu equation.
It is known that
if $v\in \mc S$ satisfies the equation \eqref{eqn_1.4}, then $v$ is an absolute minimizer for
$\mc F_A$ on $\mc S$.

\section{Classes of convex functions on $\Delta$}

In this section, we introduce several  classes of convex functions $\mc P, \mc S,\mc C_\infty$ and $\mc C_\ast$
on $\Delta$(or $\bar\Delta$).

The first class is the set of piecewise rational linear convex functions on $\bar\Delta$ which is denoted by $\mc P$.

The second class of functions is $\mc S$ that is defined in the previous section.

The third class of functions is considered by Donaldson
in \cite{D1}: let $\mc C_\infty$ be the set of
continuous convex functions on $\bar\Delta$ which are smooth in the interior.

Since we essentially consider $\mc C_\infty/\mc G$, we introduce a subset $\tilde{\mc C}\subset \mc C_\infty$ that
is isomorphic to $\mc C_\infty/\mc G$.  This can be done as the following:  fix
a point $p_o\in \Delta$ and consider the normalized functions:
$$
\tilde{\mc C}=\{u\in \mc C_\infty| u\geq u(p_o)=0\}.
$$

 For any $u\in \tilde{\mc C}$ we set
\begin{equation}
 \|u\|_b:=\int_{\partial \Delta} ud\sigma,
\end{equation}
i.e, we think it as a norm of $u$, where $b$ stands for boundary. It is known that for a sequence of functions $\{u^{(k)}\}\in
\tilde{\mc C}$
that is uniformly bounded with respect to this norm, the sequence is locally uniformly converges to a convex function $u$ in $\Delta$.
Based on this fact,
in \cite{ZZ1}, Zhou-Zhu
consider the functional $\mc F_A$ on the "compactification" of $\tilde {\mc C}$. Be precise,  they consider classes of
functions $\mc C_\ast$ and $\mc C_\ast^P$ as the following.
Let
\begin{eqnarray*}
\mc C_\ast&=&\{u| \mbox{there exist a constant $C>0$  and  a  sequence of $\{u^{(k)}\}$ in $\tilde{\mc C}$ }\\
&&\mbox{such that
   $\|u^{(k)}\|_b<C$ and
$u^{(k)}$  locally uniformly converges to} \\
&& \mbox{$u$ in $\Delta$}\}.
\end{eqnarray*}
For any $u\in \mc C_\ast,$   define $u$ on boundary as
$$u(q)=\lim_{\Delta\ni \xi\to q} u,\;\;\; q\in \partial \Delta.$$

It is known that a function $u\in \mc C_\ast$ is a convex function on $\Delta^\ast$ (the union of $\Delta$ and open codimensional-1 faces).
Let $P>0$ be a  constant, we define
$$
\mc C_\ast^P=\{u \in\mc C_\ast| \|u\|_b\leq P  \}.
$$
The functional $\mc F_A$ can be generalized to the classes $\mc C_\infty$ (cf. \cite{D1}), $\mc C_\ast$ and $\mc C^P_\ast$ (cf. \cite{ZZ1}).

Donaldson proved that
\begin{equation}
\inf_{u\in\mc C_\infty}\mc F_K(u)
=\inf_{u\in\mc S}\mc F_K(u),
\end{equation}
and
 when $P$ is large enough, Zhou-Zhu proved that (cf. Corollary 2.8 (\cite{ZZ1}))
\begin{equation}\label{eqn_1.5}
\inf_{u\in\mc C_\ast}\mc F_K(u)=\inf_{u\in\mc C_\ast^P}\mc F_K(u)
=\inf_{u\in\mc S}\mc F_K(u).
\end{equation}

The following lemma is important to relate functions between $\mc C_\infty$ and $\mc C_\ast$
with respect to the linear functional $\mc L_A$. The proof of lemma follows from the argument of \cite{D1} and  \cite{ZZ2}.
\begin{lemma}\label{lemma_1.2.12}
For any $u\in \mc C_*^P,$ there is a sequence of functions $u^{(k)}\in \mc C_\infty$
such that  $u^{(k)}$ locally uniformly converges to $u$ in $\Delta$ and
\begin{eqnarray}
&&\|u\|_b=\lim_{k\to\infty}\|u^{(k)}\|_b, \label{eqn_A}\\
&&\mc L_A(u)=\lim_{k\to\infty}\mc L_A(u^{(k)}).\label{eqn_B}
\end{eqnarray}
\end{lemma}
 { \bf Proof.} Without loss of generality, we can assume that
$0$ is the center of the mass of $\Delta $ and $u$ is normalized at the origin, i.e,
$u\geq u(0)=0$.

  Let
$$\tilde u_k(\xi)=u(r_k\xi),$$
 with  $r_k<1$ and $\lim\limits_{k\to\infty}r_k=1$.  Then  the  sequence $\{\tilde u^{(k)}\}$
in $C(\bar{\Delta})$ locally uniformly converges to u. From the construction and the convexity we know $u(\xi)\geq \tilde u^{(k)}(\xi)$
for any $\xi\in \partial \Delta.$  Then Lebesgue dominated convergence theorem implies that
$\tilde u^{(k)}(\xi)$ converges $u|_{\partial \Delta}$  in $L^1(\partial \Delta)$ as $k\to \infty.$ In particular,
\begin{equation}\label{eqn_1.14}
\lim_{k\to\infty}\int_{\partial \Delta} |\tilde u ^{(k)}-u| d\sigma =0.
\end{equation}

Let $\tilde u_\epsilon^{(k)}$ be the involution $\tilde u^{(k)}\star \rho_{\epsilon}(\xi),$ where
 $\rho_{\epsilon}\geq 0$ is a smooth mollifier of $\real^n$ whose support is in
  $B_\epsilon(0).$     For any $k $ and
  $$
  \epsilon <\frac{1-r_k}{4}dist(0,\partial \Delta), $$
  $\tilde u_\epsilon^{(k)}$ is a smooth convex function in $\bar\Delta.$
 In fact, for any $\xi,\xi'\in \bar\Delta,$ and $t\in (0,1)$
 \begin{eqnarray*}
 \tilde u_\epsilon^{(k)}(t\xi+(1-t)\xi')&=&\int_{\mathbb R^n}
 \tilde u ^{(k)}(t\xi+(1-t)\xi'-y)\rho_{\epsilon}(y)dy \\
 &=&\int_{B_\epsilon(0)}
 \tilde u^{(k)}(t(\xi-y)+(1-t)(\xi'-y))\rho_{\epsilon}(y)dy\\
 &\leq &\int_{B_\epsilon(0)}
 [t \tilde u^{(k)}(\xi-y)+(1-t)\tilde u^{(k)}(\xi'-y)]\rho_{\epsilon}(y)dy \\
 &=& t\tilde u_\epsilon^{(k)}( \xi)+(1-t)\tilde u_\epsilon^{(k)}( \xi')
 \end{eqnarray*}
 Then $\tilde u_\epsilon^{(k)}$ uniformly converges to $\tilde u^{(k)}$ on $\bar \Delta$ as $\epsilon\to 0.$
 We choose  $\epsilon_k>0$  such that
 \begin{equation}\label{eqn_1.15}
 \|\tilde u_{\epsilon_k}^{(k)}-\tilde u ^{(k)}\|_{L^\infty(\bar\Delta)}\leq \frac{1}{k}.
 \end{equation}
Let $u^{(k)}=\tilde u_{\epsilon_k}^{(k)} .$ Since
$$\int_{\partial \Delta} |u-  u ^{(k)}|d\sigma\leq
 \int_{\partial \Delta} |u-\tilde u ^{(k)}|d\sigma+\int_{\partial \Delta} |\tilde u ^{(k)}- u ^{(k)}|d\sigma,
$$
by \eqref{eqn_1.14} and \eqref{eqn_1.15}, we conclude that $u^{(k)}$ converges to $u|_{\partial \Delta}$ in $L^1(\partial \Delta).$
We prove \eqref{eqn_A}.
 On the other hand, $u^{(k)}$ locally uniformly converges to $u$ in $\Delta.$ Then \eqref{eqn_B} still holds. The lemma is proved.  q.e.d. \\

\section{K-Stability}

\subsection{Relative K-stability}\label{sect_1.2}

We start with the following definition
\begin{defn}\label{definition_1.2.1}[relatively \rm{K}-polystable]
 Let $A\in
C^{\infty}(\bar\Delta)$ be a smooth function on $\bar\Delta$.
$({\Delta},A)$ is called {\em relatively \rm{K}-polystable} if
\begin{equation}\label{eqn_4.1a}
\mc
L_{A}(u)\geq 0 \;\;\;\forall
u\in \mc P,
\end{equation}
and there is no nontrivial extremal function in $\mc P$.
\end{defn}
In \cite{ZZ2}, Zhou-Zhu proved that the relatively K-polystability is a necessary condition for the existence
of extremal metrics.

Since convex functions can be approximated by a sequence of rational piecewise-linear convex functions,
we  may assume that
\begin{equation}\label{eqn_1.6}
\mc L_A(u)\geq 0,\;\;\; \forall u\in \mc S.
\end{equation}
By Lemma \ref{lemma_1.2.12},  we also have
\begin{equation}\label{eqn_1.7}
\mc L_A(u)\geq 0,\;\;\; \forall u\in \mc C_\ast.
\end{equation}
If we replacing $\mc P$ in \eqref{eqn_4.1a} by  other classes in Definition \ref{definition_1.2.1}, such as $\mc S$, $\mc C_\infty$ or $\mc C_\ast^P$, we may define
stronger version of stability. This leads to the following definition:
\begin{defn}\label{definition_1.2.2}[relatively strong K-stability]
 Let $A\in
C^{\infty}(\bar\Delta)$ be a smooth function on $\bar\Delta$.
$({\Delta},A)$ is called {\em relatively strong \rm{K}-stable} if
\begin{equation}\label{eqn_4.2a}
\mc
L_{A}(u)\geq 0 \;\;\;\forall
u\in \mc C_\ast^P,
\end{equation}
for some large $P$,
and there is no nontrivial extremal function in $\mc C^P_\ast$.
\end{defn}



The main result of the paper is
\begin{theorem}\label{theorem_1.2.7}
If the Abreu equation \eqref{eqn_1.4} has a solution in $\mc S$, then $(\Delta,A)$ is  relative strong K-stable.
\end{theorem}

\subsection{Uniform K-stability}\label{sect_4.2}

In \cite{D1} Donaldson introduced a  new notion of stability for toric varieties which we call uniformly K-stable  here.

\begin{defn}\label{definition_1.2.4}
$({\Delta},A)$ is called uniformly K-stable if there exists a constant $\lambda>0$
such that
\begin{equation}\label{eqn_1.9}
\mc L_A(u)\geq \lambda\int_{\partial \Delta} u d \sigma,\;\;\;
\forall u\in \tilde {\mc C}.
\end{equation}
Sometimes, we say that $\Delta$ is
$(A,\lambda)$-stable.
\end{defn}

 By Lemma \ref{lemma_1.2.12}, it is clear that
 $$
 \mbox{ uniformly K-stable}
 \Rightarrow
 \mbox{relatively K-polystable}.
 $$
 The reverse is unknown in general, however there are results when $\dim \Delta=2$:
\begin{itemize}
\item when $\dim \Delta=2$ and $A\geq 0$, this is proved by Donaldson (\cite{D1});
\item when $\dim \Delta=2$ and $A$ is linear, this is proved by Wang-Zhou(\cite{WZ}).
\end{itemize}
By Lemma \ref{lemma_1.2.12}, we have
\begin{prop}
$(\Delta, A)$ is relatively strong K-stable if and only if
it is uniform K-stable.
\end{prop}
{\bf Proof. }"$\Rightarrow$". If $(\Delta, A)$ is not uniform K-stable, then for any  integer $n>0$ there exists
a function $u_n\in \tilde C$ with $\|u_n\|_b=1$ such that
$$
\mc L_A(u_n)\leq \frac{1}{n}.
$$
Then $u_n$ locally uniformly converges to a function $u\in \mc C^P_\ast$ and $\mc L_A(u)=0$. Since
$(\Delta,A)$ is  relatively strong K-stable, this implies that $u\in \mc G$. Since $u_n$ are normalized, $u\equiv 0$.
However, by the assumption, we know that
$$\int_\Delta Au_n\geq 1-\frac{1}{n}$$
and the limit of this is $\int_\Delta Au=0$, which is impossible. This completes the proof of one direction.

"$\Leftarrow$" Let $u\in \mc C^P_\ast$, then by Lemma \ref{lemma_1.2.12}, there exists a sequence of functions
$u_k\in \tilde{\mc C}$ that locally uniform converges to $u$ and $\mc L_A(u_k)\to \mc L_A(u)$. Hence $\mc L_A(u)\geq 0$.

Now suppose that $u$ is a nontrivial extremal function, i.e, $u\not\equiv 0$ and $\mc L_K(u)=0$.
Without loss of generality, by Lemma \ref{lemma_1.2.12} we assume that the limit of
$$
\lim_{k\to\infty}\|u_k\|_b\to \|u\|_b
$$
is $C>0$.
However, $\mc L_A(u)=0$ is the limit
$$
\mc L_A(u_k)\geq \lambda \|u_k\|_b\to \lambda C>0.
$$
We get a contradiction. Hence $u$ can not be a nontrivial extremal function. This completes the proof. q.e.d.

\v As a corollary of Theorem \ref{theorem_1.2.7}, we have
\begin{theorem}\label{theorem_4.6}
If the Abreu equation \eqref{eqn_1.4} has a solution in $\mc S$, then $(\Delta, A)$ is uniform K-stable.
\end{theorem}
Namely, we prove that the uniform K-stability is a necessary condition for existing a solution of \eqref{eqn_1.4} in $\mc S$.

The uniform stability should be important and it is interesting to consider its generalization to arbitrary K\"ahler manifolds. In
\cite{Sz}, the author proposed a notion of uniform stability for general K\"ahler manifold and its version on toric varieties can be read as
\begin{equation}\label{eqn_4.6a}
\mc L_A(u)\geq \lambda \|u\|,\;\;\;\forall u\in \tilde{\mc C},
\end{equation}
where $\|u\|$ is certain norm which is chosen to be $\|u\|_{L^{\frac{n}{n-1}}}$ in \cite{Sz}. From the point of view of \eqref{eqn_4.6a},
 the uniform stability
defined in this paper is to choose $\|u\|$ as $\|u\|_b$. It is shown in \cite{Sz} that the norm $\|u\|_b$ is stronger than $\|u\|_{L^{\frac{n}{n-1}}}$.

\subsection{Consequences of the uniform K-stability}\label{sect_4.3}

We mention some consequences of the uniform K-stability.

In \cite{D1}, Donaldson has a nice observation:
\begin{theorem}[Donaldson, Prop. 5.1.2\cite{D1}]\label{theorem_4.7}
If $(\Delta,A)$ is uniform stable, then the functional $\mc F_A$ is bounded from  below in $\mc C_\infty$.
\end{theorem}

By  Tian's works,  the properness of Mabuchi functional is another important notion to stress the stability. Tian conjectures
that the extremal metric exists if and only if the Mabuchi functional is proper. In \cite{ZZ1}, Zhou-Zhu
formulates the properness on toric varieties. Here, we show that the properness is a necessary condition. Be precise, we have the following result.
\begin{theorem}\label{theorem_4.8}
If the Abreu equation \eqref{eqn_1.4} has a solution $v\in \mc S$, then there exist constants $\epsilon, C>0$ such that
$$
\mc F_K(u)\geq -C+ \epsilon \int_\Delta ud\mu, \;\;\;\forall u\in \tilde{\mc C}.
$$
\end{theorem}
{\bf Proof. }Let $u_o$ be the function given by \eqref{eqn_1.0} and suppose
$$
A_o=-\sum\frac{\partial^2u_o^{ij}}{\partial x_i\partial x_j}.
$$
Then $u_o\in \mc S$ is the minimizer  of the functional $\mc F_{A_o}$ in $\mc C_\infty$, i.e,
\begin{equation}\label{eqn_4.7a}
\mc F_{A_o}(u)\geq \mc F_{A_o}(u_o)=: -C_o.
\end{equation}
Also
 there exists a constant $R$ such that
\begin{equation}
\mc L_{A_o}(u)\leq R \|u\|_b.
\end{equation}
In fact, this follows from the  fact
$$\left|\int_\Delta A_oud\mu\right|\leq |A_o|_{L^\infty}\|u\|_{L^1(\Delta)}\leq C\|u\|_b.$$

For any constant $r$ applying \eqref{eqn_4.7a} to $ru$ we have
$$
\mc F_{A_o}(ru)= -\int\log\det(u_{ij})-C(r)+r\mc L_{A_o}(u)\geq -C_o.
$$
Hence,
\begin{eqnarray*}
\mc F_{A}(u)&=& -\int\log\det(u_{ij})+\mc L_{A}(u)\\
&\geq&
-C_o+C(r)-r\mc L_{A_o}(u)+\mc L_A(u)\\
&\geq& -C_o+C(r)+
(\lambda-r R)\|u\|_b .
\end{eqnarray*}
In the last step, we use the property of uniform stability.
Choose $r$ small such that $\epsilon':=\lambda-rR>0$, we have
\begin{equation}\label{eqn_4.9a}
\mc F_A(u)\geq -C+\epsilon'\|u\|_b.
\end{equation}
By the convexity of $u$, it is straightforward to show that when $u\in \tilde{\mc C}$
there exists a constant $C'$ such that $\|u\|_b\geq C'\|u\|_{L^1(\Delta)}$. Hence,
$$
\mc F_A(u)\geq -C+\epsilon'C'\int_\Delta ud\mu,
$$
Set $\epsilon=\epsilon'C'$, we finish the proof. q.e.d.

\v
From the proof, we see that we have a stronger version of the theorem, i.e, \eqref{eqn_4.9a}.

Another consequence of Theorem \ref{theorem_4.6} is
\begin{prop}
If $v\in \mc S$ is the solution of the Abreu equation \eqref{eqn_1.4}, then there exists a constant $C$ depending on
$\lambda$
such that
$$
\|v\|_b\leq C.
$$
\end{prop}
{\bf Proof. } Note that by \eqref{eqn_5.1a}, $\mc L_A(v)=n \mbox{Area}(\Delta)=:C_1$. Then by the unform K-stability, we have
$$
\|u\|_b\leq \lambda^{-1}\mc L_A(v)=\lambda^{-1}C_1.
$$
q.e.d.

\section{Proof of Theorem \ref{theorem_1.2.7}}

Assume that $v\in \mc S$ is the solution of the Abreu equation \eqref{eqn_1.4}.

Recall that for any convex function $u$ on a domain $\Omega\subset \real^n$, it defines a Monge-Ampere measure $M_u$ on
$\Omega$.

Let $u$ be a convex  function. For any segment $I\subset\subset \Delta$, the convex function $u$ defines
a convex function $w:=u|_I$ on $I$. It defines a Monge-Ampere measure $M_w$ on $I$, we denote this by $N$.

\begin{lemma}\label{lemma_1.2.8}
Let $u\in \mc C^K_\ast$  and $u^{(k)}\in \mathcal C_\infty$ locally uniformly converges to $u.$ If $N(I)=m>0$, then $$\mc L_A(u^{(k)})> \tau m$$ for some positive constant $\tau$ independent of k.
\end{lemma}
{\bf Proof. }First assume that $u\in \mc C_\infty$.
Then
$$
\mc L_A(u)=\int_{\partial \Delta} ud\sigma -\int_\Delta A ud\mu
=\int_{\partial \Delta} ud\sigma +\int_\Delta v^{ij}_{ij} ud\mu.
$$
By Lemma 3.3.5 in \cite{D1}, (there is a misprint of sign in the paper, for the correct one, please refer to (3.3.6) in \cite{D1}),
$$
\int_\Delta v^{ij}_{ij} u=\int_{\Delta} v^{ij} u_{ij} - \int_{\partial\Delta } ud\sigma,
$$
we have
\begin{equation}\label{eqn_5.1a}
\mc L_A(u)= \int_{\Delta} v^{ij} u_{ij}d\mu.
\end{equation}

Let $p$ be the midpoint of $I$.  For simplicity, we assume that
$p$ is the origin, $I$ is on the $x_1$ axis and  $I=(-a,a)$.
 Suppose that there is a Euclidean ball $B:=B_{\epsilon_o}(0)$ in $x_1=0$ plane
such that
$I\times B\subset\subset \Delta$.  Consider the functions
$$
w_x(x_1)=u(x_1,x),\;\;\; x_1\in I, x\in B.
$$
For any $x\in B,$ $w_x$ are convex functions on $I$.   We denote the Monge-Ampere measure on $I$ induced by $w_x$ by
$N_x$. Note that $N_0=N$. By the weakly convergence of Monge-Ampere measure, we know that
there exists a small $B$ such that for any $x\in B$
\begin{equation}
N_x(I)\geq m/2.
\end{equation}
In fact, if not, then there exists a sequence of
$x_k\to 0$ such that $M_{x_k}(I)< m/2$. However, the sequence of functions
$w_{x_k}$ uniformly converges to $w_0$ in $I$, their Monge-Ampere measures should converge in any compact set, i.e,
$N_{x_k}(I)\to m$, this contradicts the assumption.

On the other hand, the eigenvalues of $v^{ij}$ are bounded below in $I\times B$, let $\delta$ be the lower bound. Then
\begin{eqnarray*}
\mc L_A(u)
 &\geq& \int_{I\times B} v^{ij}u_{ij}d\mu
\geq \delta\int_{I\times B} Trace(u_{ij})d\mu \\
&\geq& \delta\int_{I\times B} u_{11}d\mu
=\delta\int_B N_x(I)dx\geq  \frac{m\delta}{2} Vol(B).
\end{eqnarray*}
This completes the proof for $u\in \mc C_\infty$.

Now suppose that $u$ is a limit of a sequence $u^{(k)}\in \mc C_\infty$.
Then $u^{(k)}$ converges to $u$ uniformly on $I\times B$. For each, $u^{(k)}$ we
repeat the above argument: let $w^{(k)}_x$, $N_x^{(k)}$ replace $w_x$ and $N_x$. Since $w^{(k)}_x$
converges to $w_x$ uniformly with respect to $k$ and $x$, $N^{(k)}_x(I)$ also converges to
$N_x(I)$ uniformly with respect to $k$ and $x$. Hence, we conclude that when $k>K$ for some large $K$
and $x\in B$,
$$
N^{(k)}_x(I)\geq m/4.
$$
Therefore, by the same computation as above, we have
$$
\mc L_A(u^{(k)})
 \geq  \frac{m\delta}{4} Vol(B).\;\;\;\;\;\;\;q.e.d.
$$

\v \n
{\bf Proof of Theorem \ref{theorem_1.2.7}:} If  $(\Delta,A)$ is not  uniformly K-polystable, then there is a sequence of convex function $u^{(k)}$ with the property
\begin{equation}\label{eqn_1.11}
\int_{\partial \Delta} u^{(k)}d\sigma=1,\;\;\;\mbox{ and }\;\;\;
\lim_{k\to\infty} \mc L_A(u^{(k)})= 0.
\end{equation}
In particular
\begin{equation}\label{eqn_1.12}
\lim_{k\to\infty}\int_{\Delta} A u^{(k)}d\mu =1.
\end{equation}
Then $u^{(k)}$ locally uniformly  converges to a function $u\in C_*^K$.

 By the Lemma \ref{lemma_1.2.8}, we know that $\lim\limits_{k\to\infty}\mc L_A(u^{(k)})=0$ only if for any interior interval $I\subset\subset \Delta$
the Monge-Ampere measure of $u|_I$ is 0. If this is the case, $u$ must be affine linear. In fact,
suppose we normalize $u$ at some point $p$ such that $u\geq u(p)=0$. Now consider any line
$\ell$ through $p$. By the assumption, the Monge-Ampere measure of $u|_\ell$ is then trivial.
Hence $u\equiv 0 $ on $\ell$. This is true for any line, hence $u\equiv 0$.
In particular
$$\int_{\Delta}Aud\mu= 0.$$
Since $u^{(k)}$ locally uniformly converges to $u,$ by the result of Donaldson (cf. the last line in p332 of \cite{D1},) we know
\begin{equation}\label{eqn_1.13}
\lim_{k\to \infty}\int_{\Delta} Au^{(k)}d\mu =\int_{\Delta}Aud\mu .
\end{equation}
Then $\lim\limits_{k\to \infty}\int_{\Delta} Au^{(k)}d\mu=0.$ It contradicts to \eqref{eqn_1.12}.  q.e.d. \\

\section{Generalizations}\label{sect_2}
 In \cite{D2}
Donaldson also consider the case that $\Delta$ is any bounded convex domain with strictly convex smooth boundary.

  Let $\Omega \subset \mathbb R^n $ be a bounded convex domain with strictly convex smooth boundary and $\sigma$ be a smooth positive measure on $\partial \Omega$.  First we recall definition of $\mathcal S_{\Omega,\sigma}$.

 For any point $P\in\partial \Omega$ we call a coordinates
    adapted if the local coordinates
    $\xi, \eta_{1}, \dots \eta_{n-1}$ near $P$ so that $\partial \Omega$
    is given by the equation $\xi=0$ and the normal derivative of $\xi$
    on the boundary is  $\sigma^{-1}$. Define functions $\alpha_{p},p=1,2,3,$ of a positive real variable
    by
    $$ \alpha_{1}(t) = - \log t \ ,\ \alpha_{2}(t)
    = t^{-1}\ , \alpha_{3}(t) = t^{-2}. $$

    \begin{defn}
      Let $\mathcal S_{\Omega,\sigma}$ be all continuous convex functions
    $u$ on $\overline{\Omega}$ such that
    \begin{itemize}
    \item $u$ is   smooth and strictly convex on $\Omega$;
    \item In a neighbourhood of any point $p$ of $\partial \Omega$ there
    are adapted co-ordinates $(\xi,\underline{\eta_{i}})$ in which
    $$ u = \xi \log \xi + f$$
    where for $p\geq 1$, $p+q\leq 3$
    $$  \vert \nabla_{\xi}^{p} \nabla_{\underline{\eta}}^{q} f \vert=
     o(\alpha_{p}(\xi)).
    $$
    as $\xi\rightarrow 0$.
    \end{itemize}
    \end{defn}

Donaldson considers
  the problem of minimising the
  functional
    \begin{equation} {\mathcal F}(u) = - \int_{\Omega} \log \det (u_{ij})d\mu
    + \mathcal L_{\Omega,A}(u),  
    \end{equation}
where  $\Omega$ is a bounded convex domain with strictly convex smooth boundary and
    $\sigma$ is a smooth positive measure on $\partial \Omega$.
Here
$$\mathcal L_{\Omega,A}(u)=\int_{\partial
    \Omega} u d\sigma - \int_{\Omega}A ud\mu.$$

Similarly,  we can define
the relatively strong stability (cf. Definition \ref{definition_1.2.2})
and the uniform stability(cf. Definition \ref{definition_1.2.4})  for $(\Omega,A).$
Then by the same argument of proof of Theorem \ref{theorem_1.2.7} we have:
\begin{prop}
If the Abreu equation \eqref{eqn_1.4} on $\Omega$ has a solution in $\mc S_{\Omega,\sigma}$, then $(\Omega,A)$ is relatively strong stable,
hence, uniformly stable.
\end{prop}

\end{document}